\documentclass{ert-l}
\usepackage[english]{babel}
\usepackage{mathdots,enumerate,amsmath}
\usepackage{epsfig}

\issueinfo{00}{}{}{2005}
\copyrightinfo{2005}{American Mathematical Society}

\newtheorem{theorem}{Theorem}[section]
\newtheorem{proposition}[theorem]{Proposition}
\newtheorem{corollary}[theorem]{Corollary}
\newtheorem{lemma}[theorem]{Lemma}

\theoremstyle{definition}
\newtheorem{definition}[theorem]{Definition}
\newtheorem{example}[theorem]{Example}

\newcommand{\N}{\mathbb{N}}
\newcommand{\Z}{\mathbb{Z}}
\newcommand{\Q}{\mathbb{Q}}
\newcommand{\R}{\mathbb{R}}
\newcommand{\Co}{\mathbb{C}}
\newcommand{\M}{\mathbb{M}}

\newcommand{\Sl}{\mathrm{SL}(2)}
\newcommand{\Spm}{\mathrm{Sp}(2m)}
\newcommand{\Sp}{\mathrm{Sp}(2)}

\newcommand{\Id}[1]{\mathrm{Id}_{#1}}
\newcommand{\ch}{\operatorname{ch}}
\newcommand{\om}{\ensuremath{\omega}}
\newcommand{\bi}[2]{\left( \begin{array}{c} \!\!\! #1 \!\!\! \\ \!\!\! #2 \!\!\! \end{array}\right)}

\newcommand{\som}[2]{\sum\limits_{\begin{array}{c} \\*[-1.4em] \scriptstyle #1 \\*[-0.4em] \scriptstyle#2 \end{array}}}

\begin{document}

\title[Characters of the irreducible representations]{Characters of the irreducible representations with fundamental highest weight for the symplectic group in characteristic $p$}

\author{Sebastien Foulle}
\address{Institut Camille Jordan\\  Universit\'e Claude Bernard - Lyon 1\\ 69622 Villeurbanne cedex\\ France}
\email{foulle@math.univ-lyon1.fr}


\date{December 13, 2005}

\begin{abstract}
Let $K$ be an algebraically closed field of characteristic $p>0$ and let $\Spm$ be the symplectic group of rank $m$ over $K$. The main theorem of this article gives the character of the rational simple $\Spm$-modules with fundamental highest weight as an explicit alternating sum of characters of Weyl modules. One obtains several formulae for the dimensions of these simple modules, what allows us to investigate the asymptotic behavior of these dimensions, for a given $p$, when the rank $m$ is growing towards infinity. One also gets the simple Weyl modules with fundamental highest weight, and the article ends with an application to the modular representation theory of the symmetric group.
\end{abstract}

\maketitle

\section{Introduction}

First we fix some notations. Let $p$ be a prime, let $K$ be an algebraically closed field of characteristic $p$ and let $\Spm$ be the symplectic group of rank $m$ defined over $K$. We choose a maximal torus $H$ and a Borel subgroup containing $H$. Let $P=X(H)$ be the weight lattice and $P^+$ be the set of dominant weights. For $\lambda \in P^+$, we note $L(\lambda)$ the rational simple $\Spm$-module with highest weight $\lambda$ and $\Delta(\lambda)$ is the Weyl module with highest weight $\lambda$. We also note $L_m(\lambda)$, especially in section 2 and 4, when the rank $m$ takes several values or grows towards infinity. For any rational $\Spm$-module $M$, $\ch M$ denotes the character of $M$. We set $\om_0=0$ and $\om_1,\dots, \om_m$ are the fundamental weights with the usual numbering: for $1 \leq r \leq m$ one has $\ch \Delta(\om_r)=\ch \bigwedge^r V- \ch \bigwedge^{r-2} V$, where $V$ is the natural representation of $\Spm$ and $\bigwedge^k V$ its $k$th exterior power. We set $\Delta(\om_r)=\{0\}$ if $r<0$ and $\bi{n}{k}=0$ for $k<0$ or $k>n$, where $\bi{n}{k}$ denotes the binomial coefficient. Finally, for $j \in \N$, we call $p$-adic expansion of $j$ the expression $j=\sum\limits_{i \geq 0} j_i p^i$ with $j_i \in \N$ and $0 \leq j_i \leq p-1$ for all $i$.

If $G$ is a special linear group or an orthogonal group over $K$, the Weyl modules with fundamental highest weight are simple $G$-modules (see \cite{Won}). This is not the case when $G$ is a symplectic group, and the main purpose of this article is to obtain the character of the simple $\Spm$-modules $L(\om_r)$. In particular, this enables us to determine which of the Weyl modules $\Delta(\om_r)$ are simple. The following theorem gives the caracter formula for $L(\om_r)$. 

\begin{theorem}\label{theoreme} For $1 \leq r \leq m$, let $R=m+1-r=\sum\limits_{i \geq f} R_i p^{i}$ be the p-adic expansion of $R$ with $R_f \not=0$, and set $\delta=(p-R_f)p^{f}$. One has 
$$\ch L(\om_r)=\sum\limits_{j \in J} \big(\ch \Delta(\om_{r-2j})-\ch \Delta(\om_{r-2j-2\delta})\big),$$
with $J= \Big\{ j \geq 0 \ \Big\vert$ if $j=\sum\limits_{i \geq 0} j_i p^i$ is the p-adic expansion of $j$, one has $j_i=0$ when $i \leq f$ and $j_i+R_i<p$ when $i \geq f+1 \Big\}$.
\end{theorem}

\begin{example}Suppose $r=2$. One finds readily

$\ch L(\om_2)=
\left\{\begin{array}{l}
\ch \Delta(\om_2)-e^0 \mbox{ if } p \mbox{ divides } m, \\
\ch \Delta(\om_2) \mbox{ otherwise}.
\end{array}\right.$
\end{example}

We proceed as follows. Section 2 is devoted to the study of formal generating functions of dimensions of the simple $\Spm$-modules $L(\om_r)$, the main tools being dual pairs between symplectic groups and tilting modules for $\Sl$. Thanks to the main theorem of Erdmann in \cite{Erd}, we show that these series are explicit rational functions. From this we deduce a first formula for $\dim L(\om_r)$.

In section 3, independent of the preceding sections, we give two proofs of the main theorem, both of which using the decomposition matrix of the Weyl modules $\Delta(\om_r)$ found by Premet and Suprunenko (\cite{PS}). In the following section, we deduce from theorem \ref{theoreme} a first asymptotic result for the dimensions of the simple modules. We also deduce a trigonometric expression for the dimension of the simple modules $L(\om_r)$. Some simple cases of this formula, namely $p \leq 5$ and $r \geq m-p+2$, have already been obtained by Gow in \cite{Gow}. With this expression, we obtain a second asymptotic formula.

In section 5, we get the Weyl modules $\Delta(\om_r)$ which are simple, and the last section contains an application to the irreducible modular representations of the symmetric group $S_n$ : if $D^{(n-i,i)}$ is a simple $S_n$-module associated to a $p$-regular partition of $n$ in at most two parts, we express the image of $D^{(n-i,i)}$ in the Grothendieck $K_0(S_n)$ as an explicit liear combination of Specht modules. The article ends with two appendices which are used in the second proof of theorem \ref{theoreme}.

\section{Formal series}

In this section, we consider several symplectic groups at the same time, thus we use the notation $L_m(\lambda)$. We set $\dim L_{0}(\om_0)=1$ and we define the following formal generating functions.

\begin{definition}
For any $d \in \N$, we set $\chi_d(z) = \sum\limits_{n \geq 0} \dim L_{d+n}(\om_n) z^n \in \Co[[z]]$.
\end{definition}

Using dual pairs and a particular product between formal series, one can reduce the study of $\chi_d(z)$ to that of an associated series $D_d(z)$ related to tilting modules for $\Sl$. We recall some facts about dual pairs and tilting modules.

\subsection{Tilting modules}

Let $G$ be a connected simple algebraic group. A rational $G$-module $M$ is called tilting if $M$ and its dual are filtered by Weyl modules. Tilting modules have several nice properties (see \cite{Don}, \cite{Erd}, \cite{Mat}): any direct summand of a tilting module is tilting, a tensor product of tilting modules is tilting, and the indecomposable tilting modules are parametrized by the dominant weights of $G$. More precisely, each indecomposable tilting module is determined by its unique highest weight, and for any dominant weight $\lambda$, we note $T(\lambda)$ the indecomposable tilting module with highest weight $\lambda$. If $M$ is a tilting module, we denote by $[M:T(\lambda)]$ the multiplicity of $T(\lambda)$ as a direct summand in $M$.

In the following definition, we take $G=\Sl$ and $\rho=\frac{1}{2}\alpha$ where $\alpha$ is the positive root.

\begin{definition}
For $d \in \N$, we set $D_d(z)= \sum\limits_{n \geq 0} \big[T(\rho)^{\otimes n}:T(d \rho)\big] z^n \in \Co[[z]]$.
\end{definition}

The tilting modules for $\Sl$ are well-known (see \cite{Don}). This allowed Karin Erdmann \cite{Erd} to decompose the tensor powers $T(\rho)^{\otimes n}$ and to obtain the following theorem. For $k \geq 1$, $U_{k-1}$ denotes the Chebyshev polynomial of the second kind defined by  $\sin(k \theta)= U_{k-1}(\cos \theta) \sin \theta$ for $\theta \in \R$.

\begin{theorem}\cite{Erd}\label{Erdmann}
Let $d \in \N$ and let $d+1=\sum\limits_{i=f}^k d_i p^{i}$ be the $p$-adic expansion of $d+1$ with $d_f \not=0$, $d_k \not=0$. One has
$$D_d(z)=\frac{1}{z} \prod\limits_{i=f}^k \frac{U_{(p-d_i)p^{i}-1}\big(\frac{1}{2z}\big)}{U_{p^{i+1}-1}\big(\frac{1}{2z}\big)}.$$
\end{theorem}

\subsection{Dual pairs}

We recall some facts about dual pairs (see \cite{AR} and \cite{Fou} for more details). Let $A_1$ and $A_2$ be finite-dimensional $K$-algebras. A $A_1 \times A_2$-module $N$ is called a $(A_1,A_2)$ dual pair if $\mathrm{End}_{A_1}(N)=A_2$ and $\mathrm{End}_{A_2}(N)=A_1$. For any $A_1 \times A_2$-module $M$, one has $M \simeq_{A_2} \bigoplus\limits_{i \in I} S^{(i)} \otimes I^{(i)}$ with distinct indecomposable $A_2$-modules $I^{(i)}$ and vector spaces $S^{(i)}$. Thus the dimension of $S_i$ is equal to the multiplicity of $I^{(i)}$ as a direct summand in the $A_2$-module $M$. Suppose $M$ is a $(A_1,A_2)$ dual pair and set $A_1^0=\bigoplus\limits_{i \in I} \mathrm{End}{S_i}$. The algebra $A_1^0$ is semi-simple, the set of simple $A_1^0$-modules is $\{S_i\}_{i \in i}$ and $A_1^0$ is a subalgebra of $\mathrm{End}_{A_2}(M)=A_1$. It can be shown that $A_1=A_1^0 \oplus R_1$ with $R_1$ the radical of $A_1$. Thus one identifies the simple $A_1$-modules and the simple $A_1^0$-modules. This implies that the dimension of the simple $A_1$-module $S^{(i)}$ is equal to the multiplicity of the indecomposable $A_2$-module $I^{(i)}$ as a direct summand in the $A_2$-module $M$.

Let $G_1$ and $G_2$ be abstract groups and let $M$ be a $G_1 \times G_2$-module. For i=1,2, we denote by $\phi_i: G_i \rightarrow \mathrm{End}(M)$  the action of $G_i$ and by $A_i$ the subalgebra of $\mathrm{End}(M)$ generated by $\phi_i(G_i)$. The $G_1 \times G_2$-module $M$ is called a $(G_1,G_2)$ dual pair if $M$ is a $(A_1,A_2)$ dual pair.

\medskip

For $m \geq 1$, we denote by $\M$ the $(\Spm,\Sp)$ dual pair defined in \cite{AR}. The algebraic group $\Sp$ is isomorphic to $\Sl$, thus $\M$ is a $(\Spm,\Sl)$ dual pair, and it has the following property (cf. \cite{AR}): $\M$ is a tilting module for $\Sl$ isomorphic to $(\bigwedge W)^{\otimes m}$, with $W$ the natural representation of $\Sl$.

From the preceding discussion, we deduce that the dimension of simple $\Spm$-modules which are composition factors of $\M$ can be obtained by calculating the multiplicities of the indecomposable tilting $\Sl$-modules which are direct summands in $\M$.

It can be shown that $L_m(\lambda)$ is a composition factor of the $\Spm$-module $\M$ if and only if $\lambda$ is a fundamental weight or the null weight. In this case, there is a dominant weight $\tilde{\lambda}$ of $\Sl$ such that $\dim L_m(\lambda)=\big[\M:T(\tilde{\lambda})\big]$. The correspondance between $\lambda$ and $\tilde{\lambda}$ is quite simple: if $\lambda=\om_r$, then $\tilde{\lambda}=(m-r)\rho$.
Thus one has $\dim L_m(\om_r)=\big[(\bigwedge W)^{\otimes m}:T((m-r)\rho)\big]$. This result can be sharpened as follows. For $1 \leq i \leq m$ we set $\varepsilon_i=\om_i-\om_{i-1}$.

\begin{proposition}\cite{Fou} For any weight $\mu=\sum\limits_{i=1}^m \mu_i \varepsilon_i$, one has
 $$\dim L_m(\om_r)_{\mu}=\left[\bigotimes\limits^m_{i=1}\bigwedge^{1-\mu_i}W:T\big((m-r)\rho\big)\right].$$
\end{proposition}

The $\Sl$-modules $\bigwedge^0 W \simeq \bigwedge^2 W$ and $\bigwedge^1 W$ are respectively isomorphic to the trivial representation and to the natural representation of $\Sl$, therefore they are indecomposable tilting modules with respective highest weight $0$ and $\rho$, and the preceding proposition can be restated as follows.

\begin{proposition}\label{paire} One has
 $$\dim L_m(\om_r)=\left[\big(2T(0) \oplus T(\rho)\big)^{\otimes m}:T\big((m-r)\rho\big)\right]$$ 
and for $0 \leq i \leq \frac{r}{2}$, $$\dim L_m(\om_r)_{\om_{r-2i}}=\left[T(\rho)^{\otimes (m-r+2i)}:T\big((m-r)\rho\big)\right].$$
\end{proposition}

It should be noted that the character of $L_m(\om_r)$, i.e. dimensions of all weight spaces, can be obtained from proposition \ref{paire} and theorem \ref{Erdmann}, but theorem \ref{theoreme} is much more convenient and elegant.

\subsection{Binomial product}\label{convolution}

We define a certain product between formal series, which allows us in the next section to give a simple expression of $\chi_d(z)$ in function of $D_d(z)$.

\begin{definition}
For $A=\sum\limits_{n\geq 0} a_n z^n \in \Co[[z]]$ and $B=\sum\limits_{n\geq 0} b_n z^n \in \Co[[z]]$, we set $$A*B=\sum\limits_{n\geq 0} \left(\sum\limits_{k=0}^n \bi{n}{k} a_k b_{n-k}\right) z^n.$$  
\end{definition}

By Newton's binomial formula one has $\frac{1}{1-az}*\frac{1}{1-bz}=\frac{1}{1-(a+b)z}$. The following proposition is a generalization of this formula. The proof we give here was found by Wolfgang Soergel and is faster than our original proof.

\begin{proposition}\label{produit}
Let $a \in \Co$. For $S(z) \in \Co[[z]]$ one has
 $$\frac{1}{1-az} *S(z)= \frac{1}{1-az} S\left(\frac{z}{1-az}\right).$$
\end{proposition}

\begin{proof}The proposition is obviously true if $S(z)=\frac{1}{1-bz}$. By linearity, the formula holds for $S \in V$, where $V$ is the vector space spanned by the family $\left(\frac{1}{1-bz}\right)_{b \in \Co}$. Moreover $V$ is dense in $\Co[[z]]$ for the $(z)$-adic topology, as one sees by considering Vandermonde determinants. The product $*$ being continuous, this suffices to proove the claim.
\end{proof} 

\subsection{Rational expression of $\chi_d(z)$}

We give an explicit expression of $\chi_d(z)$ as a rational function, and we discuss the asymptotic behavior of $\dim L_{d+n}(\om_n)$ for given $d$ and $p$ and the rank $m$ growing towards infinity.

\begin{proposition}\label{convole}
One has
$$\chi_d(z)=\frac{1}{z^d} \left(\frac{1}{1-2z}*D_d(z)\right).$$
\end{proposition}

\begin{proof}
We consider the $\big(\mathrm{Sp}(2d+2n),\Sl\big)$ dual pair $\M$. From proposition \ref{paire}, the dimension of $L_{d+n}(\om_n)$ is the multiplicity of $T(d \rho)$ in $$\big(2T(0)\oplus T(\rho)\big)^{\otimes (d+n)}=\bigoplus_{k=0}^{d+n} \bi{d+n}{k} 2^{d+n-k} T(\rho)^{\otimes k}.$$ Thus $\dim L_{d+n}(\om_n)$ is equal to the coefficient of $z^{d+n}$ in $\frac{1}{1-2z}*D_d(z)$.
\end{proof}

Theorem \ref{Erdmann} and proposition \ref{produit} give the following corollary.

\begin{corollary}\label{fraction} Let $d \in \N$ and let $d+1=\sum\limits_{i=f}^k d_i p^i$ be the $p$-adic expansion of $d+1$. One has
 $$\chi_d(z)= \frac{1}{z^{d+1}} \prod\limits_{i=f}^k \frac{U_{(p-d_i)p^i-1}\big(\frac{1}{2z}-1\big)}{U_{p^{i+1}-1}\big(\frac{1}{2z}-1\big)}.$$
\end{corollary}

\begin{proposition}Let $p^{k}$ be the highest power of $p$ in the $p$-adic expansion of $d+1$. The poles of $\chi_d(z)$ are simple and they have the form $$\frac{1}{2+2 \cos\left(\frac{j\pi}{p^{k+1}}\right)}=\frac{1}{4\cos^2\left(\frac{j\pi}{2p^{k+1}}\right)}, \ 1 \leq j \leq p^{k+1}-1.$$ The smallest pole is $\frac{1}{4\cos^2\left(\frac{\pi}{2p^{k+1}}\right)}$.
\end{proposition}

\begin{proof}
One has $$\chi_d(z)=\frac{1}{z^{d+1}} \frac{U_{(p-d_f)p^f-1}\left(\frac{1}{2z}-1\right)}{U_{p^{k+1}-1}\left(\frac{1}{2z}-1\right)}\prod\limits_{i=f}^{k-1} \frac{U_{(p-d_{i+1})p^{i+1}-1}\left(\frac{1}{2z}-1\right)}{U_{p^{i+1}-1}\left(\frac{1}{2z}-1\right)}$$ and for any positive integers $r$ and $s$ such that $r$ divides $s$, in particular for $r=p^{i+1}$ and $s=(p-d_{i+1})p^{i+1}$, $U_{r-1}(X)$ divides $U_{s-1}(X)$: if $s=qr$ one has $U_{s-1}(X)=U_{q-1}\big(T_r(X)\big) U_{r-1}(X)$, where $T_r(X)$ is the Chebyshev polynomial of the first kind defined by $\cos(r \theta)= T_{r}(\cos \theta)$ for $\theta \in \R$. Hence the rational fraction $\frac{U_{(p-d_{i+1})p^{i+1}-1}\left(\frac{1}{2z}-1\right)}{U_{p^{i+1}-1}\left(\frac{1}{2z}-1\right)}$ is a polynomial in $\frac{1}{z}$. Besides the $n-1$ distinct roots of $U_{n-1}(X)$ are the numbers $\cos\left(\frac{l \pi}{n}\right)$, $1 \leq l \leq n-1$, and the rational function $\chi_d(z)$ is a formal series, thus $0$ is not a pole. This establishes the first part of the proposition.

Let $z_0=\frac{1}{4\cos^2\left(\frac{\pi}{2p^{k+1}}\right)}$. One has immediately $\prod\limits_{i=f}^k U_{(p-d_i)p^i-1}\left(\frac{1}{2z_0}-1\right) \not=0$ and $\prod\limits_{i=f}^k U_{p^{i+1}-1}\left(\frac{1}{2z_0}-1\right)=0$, thus $z_0$ is a pole and it is the smallest one.\end{proof}

For any $a \in \Co^*$, one has $\frac{1}{z-a}=\frac{-1}{a}\frac{1}{1-\frac{z}{a}}=\frac{-1}{a}\sum \frac{1}{a^n}z^n$. Thus when $n \rightarrow +\infty$ the coefficient of $z^n$ in $\chi_d(z)$ is equivalent to $\frac{c}{z_0^n}$ where $z_0$ is the smallest pole and $c$ depends on $p$ and $d$. Hence we have $\dim L_{d+n}(\om_n) \mathop{\sim}\limits_{n \rightarrow \infty} c \ 4^n \cos^{2n}\left(\frac{\pi}{2p^{k+1}}\right)$. The exact value of $c$ will be obtained in proposition \ref{asymptotique}.

\subsection{First dimension formula}

We obtain an explicit formula for the dimensions of the simple $\Spm$-modules $L(\om_r)$, which provides a simple and fast algorithm for the determination of the dimensions of the simple modules $L(\om_r)$.

\begin{definition}\label{notation} For any integer $r$ with $0 \leq r \leq m$, we set $R=m+1-r$ and $R=\sum\limits_{i=f}^k R_i p^i$ is the $p$-adic expansion of $R$ with $R_f \not=0$, $R_k \not=0$.
\end{definition}

\begin{proposition}\label{serieentiere} The dimension of $L(\om_r)$ is the coefficient of $X^r$ in the formal series
 $(1-X)(1+X)^{2m+1}\prod\limits_{i=f}^k \frac{X^{2(p-R_i)p^i}-1}{X^{2p^{i+1}}-1}$.
\end{proposition}

\begin{proof} We set $d=m-r$ and we follow the method of \cite{Erd}. One has to show that the coefficient of $z^r$ in $\chi_d(z)$
is equal to the coefficient of $X^r$ in the formal series $\Theta_{d,r}(X)=(1-X)(1+X)^{2d+2r+1}\prod\limits_{i=f}^k \frac{X^{2(p-d_i)p^{i}}-1}{X^{2p^{i+1}}-1}$. Now the coefficient of $z^r$ in $\chi_d(z)$ is equal to $\frac{1}{2 i \pi} \oint \frac{\chi_d(z)}{z^{r+1}} dz$ (with a small contour of index 1 around zero).

 Let $\theta \in \R$ and define $w$ and $z$ by $w=e^{i \theta}$ and $\frac{1}{2z}-1=\cos \theta=\frac{1}{2} \left(w+\frac{1}{w}\right)$. One has $z=\frac{w}{(w+1)^2}$ and $ \sin(k \theta)= \frac{w^k-w^{-k}}{2i}$ for any $k \in \N$. From corollary \ref{fraction}, we get \begin{align*}\chi_d(z) & =\frac{(w+1)^{2d+2}}{w^{d+1}} \prod\limits_{i=f}^k \frac{\sin\big((p-d_i)p^i\theta\big)}{\sin(p^{i+1}\theta)}
\!=\!\frac{(1+w)^{2d+2}}{w^{d+1}} \prod\limits_{i=f}^k \frac{w^{(p-d_i)p^i}-w^{-(p-d_i)p^{i}}}{w^{p^{i+1}}-w^{-p^{i+1}}}\\
& =\frac{(1+w)^{2d+2}}{w^{d+1}} \prod\limits_{i=f}^k \frac{w^{2(p-d_i)p^i}-1}{w^{2p^{i+1}}-1} w^{d+1} =(1+w)^{2d+2} \prod\limits_{i=f}^k \frac{w^{2(p-d_i)p^i}-1}{w^{2p^{i+1}}-1}.
\end{align*}

Since $dz=\frac{1-w}{(1+w)^3}dw$, the map $w \mapsto z$ is biholomorphic in a small neighbourhood of zero, and it maps zero to zero. Thus the change of variable $z=\frac{w}{(w+1)^2}$ gives

 $$\frac{1}{2 i \pi} \oint \frac{\chi_d(z)}{z^{r+1}} dz
= \frac{1}{2 i \pi} \oint  \frac{(1-w)(1+w)^{2r-1}}{w^{r+1}} \chi_d\left(\frac{w}{(w+1)^2}\right) \, dw$$

$$= \frac{1}{2 i \pi} \oint  \frac{(1-w)(1+w)^{2d+2r+1}}{w^{r+1}}  \,\prod\limits_{i=f}^k \frac{w^{2(p-d_i)p^i}-1}{w^{2p^{i+1}}-1} dw= \frac{1}{2 i \pi} \oint 
 \frac{\Theta_{d,r}(w)}{w^{r+1}} dw.$$
\end{proof}

\begin{example}
Let $m=10$ and $p=2$. We obtain the following array for the dimensions of the simple modules and Weyl modules with fundamental highest weight.

\medskip

\noindent
$\begin{array}{c|c|c|c|c|c|c|c|c|c|c}
r & 1 & 2 & 3 & 4 & 5 & 6 & 7 & 8 & 9 & 10 \\
\hline
\!\! \dim L(\om_r) & 20 & 188 & 1120 & 4466 & 14344 & 29448 & 62016 & 53296 & 76096 & 1024 \!\!\!\! \\
\hline
\!\! \dim \Delta(\om_r) & 20 & 189 & 1120 & 4655 & 14364 & 33915 & 62016 & 87210 & 90440 & 58786 \!\!\!\!
\end{array}
$

\medskip

We note that $\dim L(\om_m)=1024=2^{10}$, what is a particular case of example \ref{exemple}. We also note that we have the $p$-adic expansion $m+1=1+2+2^3$, and the values of $r$ for which $\Delta(\om_r)$ is simple are $1,1+2,1+2+2^2$, as stated in corollary \ref{irredp2}.
\end{example}

\section{Two Proofs of theorem \ref{theoreme}}

We show here that theorem \ref{theoreme} is equivalent to the theorem of Premet and Suprunenko \cite{PS} which describes the composition factors of $\Delta(\om_r)$. For this purpose, it suffices to check that two particular matrices are inverse of each other. The first proof relies on a completely similar situation for modular representations of $\Sl$, while the second proof uses the fractal patterns of these matrices. Although this second proof could be omitted, we believe that the description of these fractal patterns is interessing in its own. The precise study of these patterns, which is rather technical, is postponed in the appendices. 

\subsection{Notations}

\begin{definition}Let $a \geq 0$ and $b>0$ be integers with $p$-adic expansions $a=\sum a_i p^i$, $b=\sum\limits_{i \geq s} b_i p^i$ and $b_s \not=0$. We write

\begin{enumerate}
\item $ a \subset b$ if one has $a_i =0$ or $a_i=b_i$ for any $i$.
\item $a \prec_1 b$ if $a_0= \dots = a_s=0$ and $a_i+b_i<p$ for any $i>s$.
\item $a \prec_{-1} b$ if $a_0= \dots = a_{s-1}=0$, $a_s+b_s=p$ and $a_i+b_i<p$ for any $i>s$.
\item  $a \prec b$ if $a \prec_1 b$ or $a \prec_{-1} b$.
\end{enumerate}
\end{definition}

Note that $0 \subset b$ and $0 \prec_1 b$ for any $b > 0$, whereas $0 \not\prec_{-1} b$ for any $b>0$.

The following theorem gives the decomposition matrix of the Weyl modules $\Delta(\om_r)$ for $\Spm$. It was obtained for $p>2$ by Premet and Suprunenko \cite{PS} and for any prime $p$ by Baranov and Suprunenko \cite{BS}. In the statement of the theorem, $j$ denotes a integer such that $0 \leq j \leq r$ and $r-j$ is even.

\begin{theorem}\label{composition} One has $\ch \Delta(\om_r)=\sum\limits_{\frac{r-j}{2}\subset m+1-j} \ch L(\om_j).$
\end{theorem}

 We restate theorem \ref{theoreme} as follows.

\medskip

\noindent {\bf Theorem 1.1}
{\it One has $\ch L(\om_r)=\sum\limits_{\frac{r-j}{2}\prec_1 m+1-r} \ch \Delta(\om_j)-\sum\limits_{\frac{r-j}{2}\prec_{-1} m+1-r} \ch \Delta(\om_j).$}

\medskip

These theorems lead us to the following definition.

\begin{definition}
For any integer $n \geq 1$, we define the lower unipotent triangular matrices $A(n)_{1 \leq k,l \leq n}$ and $B(n)_{1 \leq k,l \leq n}$ by:
\begin{enumerate} 
\item $A(n)_{k,l}=1$ if $\frac{k-l}{2} \prec_1 n+1-k$,
\item $A(n)_{k,l}=-1$ if $\frac{k-l}{2} \prec_{-1} n+1-k$,
\item $A(n)_{k,l}=0$ otherwise,
\item $B(n)_{k,l}=1$ if $\frac{k-l}{2} \subset n+1-l$,
\item $B(n)_{k,l}=0$ otherwise. 
\end{enumerate}
\end{definition}

Note that $A(n)_{k,l}=B(n)_{k,l}=0$ if $k-l$ is odd.

\medskip

 Assume that $A(m+1)$ and $B(m+1)$ are inverse of each other. By setting  $r=k-1$ and $j=l-1$ in the two preceding theorems, one sees that they are equivalent. It will be convenient to consider the matrices $\tilde{A}(n)$ and $\tilde{B}(n)$ defined by $\tilde{A}(n)_{u,v}=A(n)_{n+1-u,n+1-v}$ and $\tilde{B}(n)_{u,v}=B(n)_{n+1-u,n+1-v}$, that is
\begin{enumerate} 
\item $\tilde{A}(n)_{u,v}=1$ if $\frac{v-u}{2} \prec_1 u$,
\item $\tilde{A}(n)_{u,v}=-1$ if $\frac{v-u}{2} \prec_{-1} u$,
\item $\tilde{A}(n)_{u,v}=0$ otherwise,
\item $\tilde{B}(n)_{u,v}=1$ if $\frac{v-u}{2} \subset v$,
\item $\tilde{B}(n)_{u,v}=0$ otherwise. 
\end{enumerate}

Of course, the matrices $A(n)$ and $B(n)$ are inverse of each other if and only if $\tilde{A}(n)$ and $\tilde{B}(n)$ are inverse of each other. Clearly, one has

$ \tilde{A}(n+1)=\left(\begin{array}{ccc|c}
 & & &  *  \\
 & \tilde{A}(n) &  & \vdots \\
 & & &  * \\
\hline
0 & \cdots & 0 & 1 
\end{array}\right)$,
$ \tilde{B}(n+1)=\left(\begin{array}{ccc|c}
 & & &  *  \\
 & \tilde{B}(n) &  & \vdots \\
 & & &  * \\
\hline
0 & \cdots & 0 & 1 
\end{array}\right)$ and thus
$ \tilde{A}(n+1)\tilde{B}(n+1) =\left(\begin{array}{ccc|c}
 & & &  *  \\
 & \tilde{A}(n)\tilde{B}(n) &  & \vdots \\
 & & &  * \\
\hline
0 & \cdots & 0 & 1 
\end{array}\right)$.

In order to prove that $\tilde{A}(n)$ and $\tilde{B}(n)$ are inverse of each other, it is enough to show that $\tilde{A}(n+1)$ and $\tilde{B}(n+1)$ are inverse of each other. Thus theorem~\ref{theoreme} will be proved if one shows that $\tilde{A}(p^n)^{-1}=\tilde{B}(p^n)$ for any $n$ (or equivalently $A(p^n)^{-1}=B(p^n)$ for any $n$) , or if $\tilde{A}(p^n-1)^{-1}=\tilde{B}(p^n-1)$ for any $n$.

\subsection{$\Sl$-theorical proof}

We express the character of the simple $\mathrm{SL}(2)$-modules as an explicit alternating sum of characters of Weyl modules. Although the characters of the simple $\mathrm{SL}(2)$-modules are well-known, the expression we obtain seems to be new. From this we deduce that $\tilde{A}(p^n-1)^{-1}=\tilde{B}(p^n-1)$ for any $n$.

The composition factors of Weyl modules for $\Sl$ have been obtained by Carter and Cline \cite{CC}. We first give a description of these composition factors that uses the relation $\subset$. Let us recall the following result of Winter.

\begin{theorem}\cite{Win} Let $r,s \geq 0$ and let $s=\sum\limits s_i p^i$ be the $p$-adic expansion of $s$. One has $[\Delta(r \rho):L(s \rho)]=1$ if and only if there exists a finite subset $I$ of $\N$ such that $r=\sum\limits_{i \in I} (2p-2-s_i)p^i+\sum\limits_{i \not\in I} s_i p^i$.
\end{theorem}

\begin{corollary} Let $0 \leq s \leq r \leq p^k-1$. One has $[\Delta(r \rho):L(s \rho)]=1$ if and only if $\frac{r-s}{2} \subset p^k-1-s$.
\end{corollary}

\begin{proof} Set $z=\frac{r-s}{2}$ and $s=\sum\limits_{i=0}^{n-1} s_i p^i$ with $0 \leq s_i \leq p-1$ for any $i$. If $z \subset p^n-1-s$, there exists a set $I \subset \{0,1, \dots ,n-1\}$ such that $z=\sum\limits_{i \in I}(p-1-s_i)p^i$, and one has $r=s+2z=\sum\limits_{i \not\in I}s_ip^i+\sum\limits_{i \in I}(2p-2-s_i)p^i$, whence $[\Delta(r \rho):L(s \rho)]=1$ from the theorem. Reciprocally, if there exists a finite subset $I$ of $\N$ with $r=\sum\limits_{i \in I} (2p-2-s_i)p^i+\sum\limits_{i \not\in I} s_i p^i$, one has $z=\sum\limits_{i \in I} (p-1-s_i)p^i\leq p^n-1$, which implies $I \subset \{0,1, \dots ,n-1\}$ and thus $z \subset p^n-1-s$.
\end{proof}

Note that in \cite{BS}, Baranov and Suprunenko have studied the links between the combinatorics of Carter and Cline \cite{CC} and the relation $\subset$, and our preceding corollary is also easily deduced from proposition 2.8 of \cite{BS}.

For $1 \leq s' \leq p^n-1$ and $1 \leq r' \leq p^n-1$, the corollary shows that $L((s'-1) \rho)$ is a composition factor of $\Delta((r'-1) \rho)$ if and only if $\frac{r'-s'}{2} \subset p^n-s'$. This gives the following corollary.

\begin{corollary}
The decomposition matrix of the Weyl modules $\Delta(r \rho)$ with $0 \leq r \leq p^n-2$ is $B(p^n-1)$.
\end{corollary}

Recall the following result (formula (7) in paragraph I.2.16 of \cite{Jan}).

\begin{proposition}For $r \in \N$ and $0 \leq k \leq r$, $(r-2k)\rho$ is a weight of $L(r \rho)$ if and only if $p$ does not divide $\bi{r}{k}$.
\end{proposition}

We can now state the analog of theorem \ref{theoreme} for $\Sl$.

\begin{corollary}
For  $r \geq 0$, let $r+1=\sum\limits_{i \geq s}a_ip^i$ be the $p$-adic expansion of $r+1$ with $a_s \not=0$ and set $J=\left.\left\{\sum\limits_{i > s}c_ip^i \ \right| \forall i, \ 0 \leq c_i \leq a_i \right\}$. One has

$$\ch L(r \rho)=\!\!\!\!\!\!\som{j \in J}{r-2j \geq 0}\!\!\!\!\!\! \ch \Delta\big((r-2j)\rho\big)-\!\!\!\!\!\!\!\!\!\!\!\!\som{j \in J}{r-2a_sp^s-2j \geq 0}\!\!\!\!\!\!\!\!\!\!\!\! \ch \Delta\big((r-2a_sp^s-2j)\rho\big).$$ 
\end{corollary}

\begin{proof}
It is well-known (cf. \cite{Jam}, lemma 22.4) that $p$ does not divide $\bi{r}{k}$ if and only if $r_i \geq k_i$ for all $i$, where $r=\sum r_ip^i$ and $k=\sum k_ip^i$ are the $p$-adic expansions of $r$ and $k$. One has $r_0=r_1=\dots=r_{s-1}=p-1$, $r_s=a_s-1$ and $r_i=a_i$ for any $i>s$, hence by setting $x=\sum\limits_{0 \leq i \leq s}k_ip^i$, one gets that  $(r-2k)\rho$ is a weight of $L(r \rho)$ if and only if $k=x+j$ with $0 \leq x < a_s p^s$ and $j \in J$. Without lose of generality, we can assume that $r-2k \geq 0$, and it is clear that $(r-2k)\rho$ is a weight occuring in $\ch \Delta\big((r-2j)\rho\big) - \ch \Delta\big((r-2a_sp^s-2j)\rho\big)$ if and only if $k=x+j$ with $0 \leq x < a_s p^s$.
\end{proof}

\begin{corollary}For any positive integer $n$, the matrices $B(p^n-1)$ and $A(p^n-1)$ are inverse of each other.
\end{corollary}

\begin{proof} It suffices to check that if $1 \leq r' \leq p^n-1$, one has $\ch L\big((r'-1)\rho\big)= \sum\limits_{1 \leq s' \leq r'} A(p^n-1)_{r',s'} \ch \Delta\big((s'-1)\rho\big)$. Let $r'=\sum\limits_{i \geq s}a_ip^i$ be the $p$-adic expansion of $r'$ with $a_s \not=0$ and set $J=\left.\left\{\sum\limits_{i > s}c_ip^i \ \right| \forall i, \ 0 \leq c_i \leq a_i \right\}$. The preceding corollary gives

$$\ch L((r'-1)\rho)=\!\!\!\!\!\!\som{j \in J}{r'-1-2j \geq 0}\!\!\!\!\!\! \ch \Delta\big((r'-1-2j)\rho\big)-\!\!\!\!\!\!\!\!\!\!\!\!\som{j \in J}{r'-1-2a_sp^s-2j \geq 0}\!\!\!\!\!\!\!\!\!\!\!\! \ch \Delta\big((r'-1-2a_sp^s-2j)\rho\big),$$ 
and one has $(p^n-1)+1-r'=(p-a_s)p^s+\sum\limits_{i=s+1}^{n-1}(p-1-a_i)p^i$.
Hence if $1 \leq s' \leq r'$,  we have $\frac{r'-s'}{2} \prec_1 (p^n-1)+1-r'$ if and only if $s'-1=r'-1-2j$ with $j \in J$, and $\frac{r'-s'}{2} \prec_{-1} (p^n-1)+1-r'$ if and only if $s'-1=r'-1-2a_sp^s-2j$ with $j \in J$.
\end{proof}

\subsection{Matricial proof}

In this subsection, we explain how to construct the matrices $A\big(p^{n+1}\big)$ and $B\big(p^{n+1}\big)$ by using $A(p^n)$ and $B(p^n)$, and we show why this construction implies $A(p^n)^{-1}=B(p^n)$ for any $n$. The proofs are postponed in the appendices. For $n \geq 1$, we denote by $\Id{n} \in \mathcal{M}_n(K)$ the identity matrix and we set 
$ E_n=
\left(\begin{array}{ccccc}
 0 & \cdots & \cdots & \cdots & 0 \\ 
 \vdots & & & \iddots & 1 \\
 \vdots& & \iddots &  \iddots & 0 \\
\vdots & \iddots & \iddots & \iddots & \vdots \\
0 & 1 & 0 & \cdots & 0 
\end{array}\right) \in \mathcal{M}_n(K).$ 

We claim that:
\begin{enumerate}
\item $A(p)=B(p)=\Id{p}$,
\item $ A(p^{n+1})=
\left(\begin{array}{cccccc}
 A(p^n) & 0 & \cdots & \cdots & \cdots & 0 \\
 -A(p^n)E_{p^n} & \ddots & \ddots & & & \vdots \\
 A(p^n)E_{p^n}^2 & \ddots & & \ddots & & \vdots \\
 -A(p^n)E_{p^n} & \ddots & & & \ddots & \vdots \\
 A(p^n)E_{p^n}^2 & \ddots & & & & 0 \\
\vdots & & & \cdots & -A(p^n)E_{p^n}  & A(p^n) \\
\end{array}\right),$

\item $ B(p^{n+1})=
\left(\begin{array}{cccccc}
 B(p^n) & 0 & \cdots &\cdots &\cdots & 0 \\
 E_{p^n}B(p^n) & \ddots & \ddots & & & \vdots \\
 0 & \ddots & & \ddots & & \vdots \\
 \vdots & \ddots & & & \ddots & \vdots \\
 \vdots & & \ddots & & & 0 \\
 0 & \cdots & \cdots & 0 & E_{p^n}B(p^n) & B(p^n)
\end{array}\right).$
\end{enumerate}

\begin{example}
For $p=3$ and $n=4$, the matrices $A(p^n)$ and $B(p^n)$ are shown on pages \pageref{matriceinv} and \pageref{matricedecomp}. The null coefficients have not been printed. 

\begin{figure*}\label{matriceinv}
 \centering
 \epsfig{figure=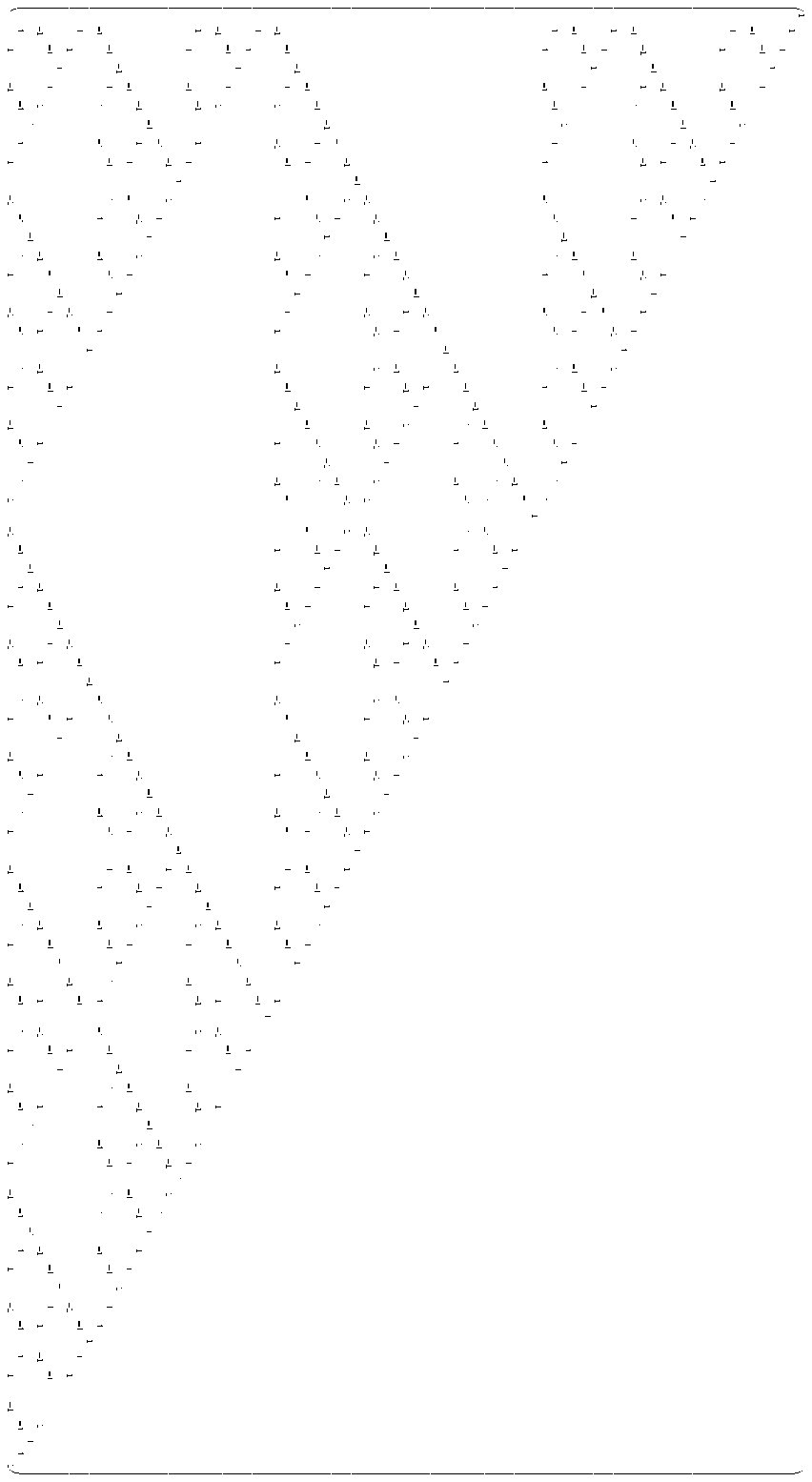,width=.8\textwidth}
\end{figure*}

\begin{figure*}\label{matricedecomp}
\centering
\epsfig{figure=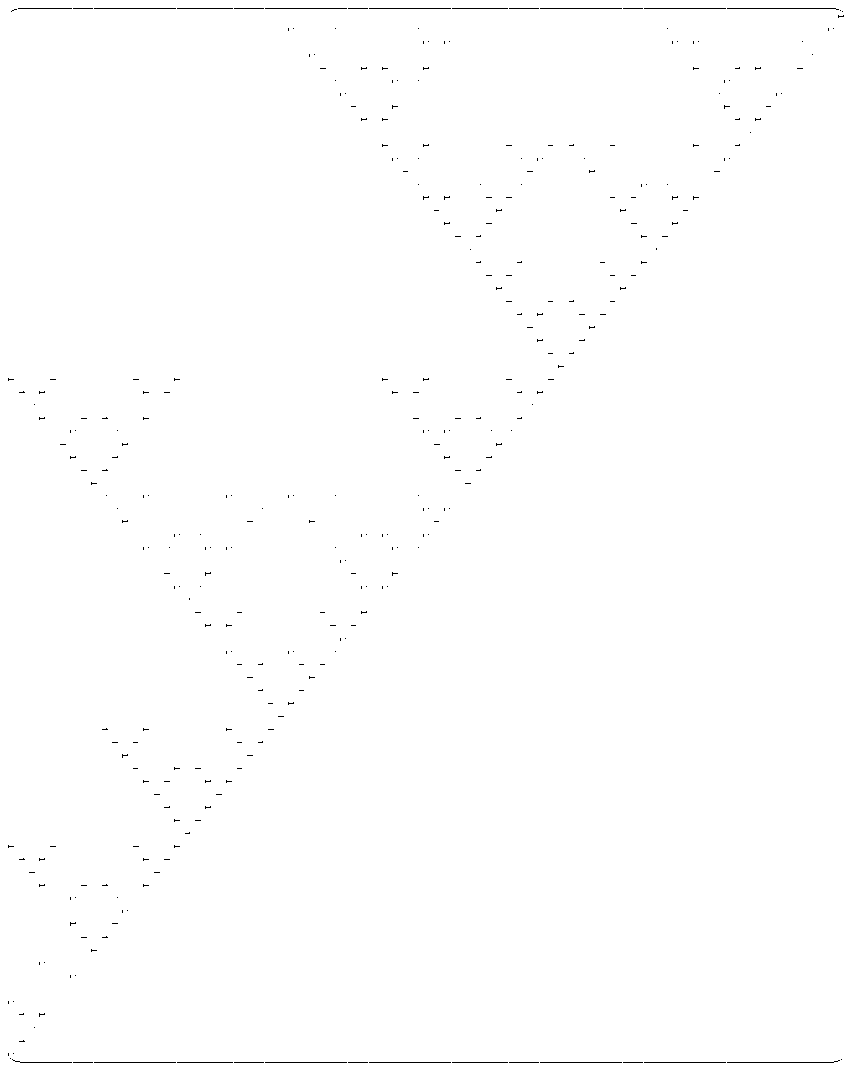,width=\textwidth}
\end{figure*}
\end{example}

From this iterative description, one obtains $A(p^n)B(p^n)=\Id{p^n}$ for any $n$: let $A$ be a noncommutative ring and let $x,y,z \in A$. One has

\medskip

\noindent
${\scriptscriptstyle \left(\!\!\!\!\begin{array}{ccccc}
x & 0 & \cdots & \cdots & 0 \\
-xz & \ddots & \ddots & & \vdots \\
xz^2 & \ddots & & \ddots & \vdots \\
-xz^3 & & & & 0 \\
\vdots & & \cdots & -xz & x
\end{array}\!\!\right)
\left(\!\!\begin{array}{ccccc}
y & 0 & \cdots & \cdots & 0 \\
zy & \ddots & \ddots & & \vdots \\
0 & \ddots & & \ddots & \vdots \\
\vdots & \ddots & \ddots & \ddots & 0 \\
0 & \cdots & 0 & zy & y
\end{array}\!\!\right)=
\left(\!\!\begin{array}{ccccc}
xy & 0 & \cdots & \cdots & 0 \\
0 & \ddots & \ddots & & \vdots \\
\vdots & \ddots & & \ddots & \vdots \\
\vdots & & \ddots & \ddots & 0 \\
0 & \cdots & \cdots & 0 & xy
\end{array}\!\!\right)}.$

\medskip
Setting $x=\!A(p^n)$, $y=B(p^n)$ and $z=E_{p^n}=z^3$, one obtains $A\big(p^{n+1}\big)B\big(p^{n+1}\big)=\Id{p^{n+1}}$ if and only if
$A(p^n)B(p^n)=\Id{p^n}$, and this is obviously true for $n=1$.

\newpage 

\section{Trigonometric formula for dimensions}

In this section, we use the notations of definition \ref{notation}. First we get the asymptotic behavior of $\dim L_m(\om_r)$ for given $p$ and $r$ when $m \rightarrow +\infty$. Then we obtain a first dimension formula which makes use of periodic sums of binomial coefficients, from which one deduces a trigonometric expression for the dimension of the simple $\Spm$-module $L(\om_r)$. It should be noticed that this second formula easily gives the partial fraction decomposition of $\chi_d(z)$. Finally we get the asymptotic behavior of $\dim L_{d+n}(\om_n)$ for given $p$ and $d$ when $n \rightarrow +\infty$.

\begin{proposition} For given $r$ and $p$, one has  $$\dim L_m(\om_r) \mathop{\sim}\limits_{m \rightarrow \infty} \dim \Delta_m(\om_r) \mathop{\sim}\limits_{m \rightarrow \infty} \frac{2^r}{r!}m^r.$$
\end{proposition}

\begin{proof}With the notations of theorem \ref{theoreme}, one has
$$\dim L(\om_r)=\sum\limits_{j \in J} \big(\dim \Delta(\om_{r-2j})-\dim \Delta(\om_{r-2j-2\delta})\big),$$
the dimension of  each term of the sum is a polynomial in $m$, and $J$ has at most $r+1$ elements. Moreover one has $\dim \Delta(\om_k)=\bi{2m}{k}-\bi{2m}{k-2} \mathop{\sim}\limits_{m \rightarrow \infty} \frac{2^k}{k!}m^k$ for any $k$, and the term of highest degree is $\Delta(\om_r)$.
\end{proof}

\begin{definition}\label{notationbis}
Set $A=\left\{\sum\limits_{i=f+1}^k a_ip^i \ \Bigg| \ \forall i, \ 0 \leq a_i \leq p-1-R_i \right\}$ if $f<k$, and $A=\{0\}$ if $f=k$.
\end{definition}

\begin{proposition}\label{sommebinome}
One has $$\dim L(\om_r)=\sum\limits_{a \in A}\sum\limits_{n \in \Z} \left(\bi{2m}{r-2a+2np^{k+1}}-\bi{2m}{r-2-2a+2np^{k+1}}\right).$$
\end{proposition}

\begin{proof}
From theorem \ref{theoreme}, the dimension of $L(\om_r)$ is equal to

 \begin{gather*}\sum\limits_{j \in J}\left(\bi{2m}{r-2j}+\bi{2m}{r-2j-2\delta-2}-\bi{2m}{r-2j-2}-\bi{2m}{r-2j-2\delta}\right)\\
=\sum\limits_{j \in J}\scriptstyle\left(\bi{2m}{r-2j}+\bi{2m}{2m-r+2j+2\delta+2}-\bi{2m}{r-2j-2}-\bi{2m}{2m-r+2j+2\delta}\right).\end{gather*}

\noindent It is enough to show that 

$$\sum\limits_{j \in J}\left(\bi{2m}{r-2j}+\bi{2m}{2m-r+2j+2\delta+2}\right)=\som{a \in A}{ n \in \Z}\bi{2m}{r-2a+2np^{k+1}},$$
 the second identity 
$$\sum\limits_{j \in J}\left(\bi{2m}{r-2j-2}+\bi{2m}{2m-r+2j+2\delta}\right)=\som{a \in A}{ n \in \Z}\bi{2m}{r-2-2a+2np^{k+1}}$$ 
being obtained in a similar manner.

 Set $S_1=\sum\limits_{j \in J}\bi{2m}{r-2j}$ and  $S_2=\sum\limits_{j \in J}\bi{2m}{2m-r+2j+2\delta+2}$. 
One has $J=\big\{a+p^{k+1}n \ \big| \ a \in A, n \in \N \big\}$, thus $S_1=\som{a \in A}{ n \in \N}
\bi{2m}{r-2a-2np^{k+1}}$. 
Moreover, the mapping which sends $a \in A$ onto $\left(p^{k+1}-p^{f+1}-\sum\limits_{i=f+1}^k R_ip^i\right)-a$ is an involution of $A$: for $a=\sum\limits_{i=f+1}^k a_ip^i$ with  $0 \leq a_i \leq p-1-R_i$ when $ f+1 \leq i \leq k$, one has $ \left(p^{k+1}-p^{f+1}-\sum\limits_{i=f+1}^k R_ip^i\right)-a=\sum\limits_{i=f+1}^k(p-1-R_i-a_i)p^i$ which is an element of~$A$.

\noindent Thus $S_2=\!\!\!\!\som{a \in A}{ n \in \N}\!\!\!
\bi{2m}{2m-r+2\delta+2+2p^{k+1}n+2\left(p^{k+1}-p^{f+1}-\sum\limits_{i=f+1}^k R_ip^i-a\right)}$. 

Now 
\begin{align*}
2m&-r+2\delta+2+2p^{k+1}n+2\left(p^{k+1}-p^{f+1}-\sum\limits_{i=f+1}^k R_ip^i-a\right)\\
&=r+2(m-r+1)+2\delta+2p^{k+1}n +2p^{k+1}-2p^{f+1}-2\sum\limits_{i=f+1}^k R_ip^i-2a\\
&=r+2R+2p^{k+1}n+2(p-R_f)p^f+2p^{k+1}-2p^{f+1}-2\sum\limits_{i=f+1}^k R_ip^i-2a\\
&=r+2p^{k+1}-2a+2p^{k+1}n=r+2p^{k+1}(n+1)-2a.\end{align*}

Hence $S_2=\som{a \in A}{ n \in \N}\bi{2m}{r-2a+2(n+1)p^{k+1}}$, and the proposition is proved.
\end{proof}

The following technical lemma allows us to compute the sums of binomials coefficients that appear in the proposition.

\begin{lemma}\label{lemmetechnique}
Let $q,s \in \N$ and $r \in \Z$. One has
$$\sum\limits_{n \in \Z}\bi{q}{r+ns}=\frac{1}{s}\sum\limits_{j=1}^s \cos\left(\frac{j \pi(q-2r)}{s}\right)\left(2 \cos\left(\frac{j \pi}{s}\right)\right)^q.$$
\end{lemma}

\begin{proof}
We set $\eta=e^{\frac{i\pi}{s}}$ and we first show that 
$\sum\limits_{n \in \Z}\bi{q}{r+ns}=\frac{1}{s}\sum\limits_{j=1}^s \eta^{-2rj} \big(1+~\eta^{2j}\big)^q$.

 One has \begin{align*}\sum\limits_{j=1}^s \eta^{-2rj} (1+\eta^{2j})^q =\sum\limits_{j=1}^s \eta^{-2rj}\sum\limits_{l \in \Z} \bi{q}{l} \eta^{2jl}&=\sum\limits_{l \in \Z} \bi{q}{l} \sum\limits_{j=1}^s \eta^{2j(l-r)}\\
= \!\!\!\! \sum\limits_{l \equiv r \!\! \mod s }\!\!\bi{q}{l} \sum\limits_{j=1}^s \eta^{2j(l-r)}+\!\!\!\!\sum\limits_{l \not\equiv r \!\! \mod s }\!\!\bi{q}{l} \sum\limits_{j=1}^s \eta^{2j(l-r)}&= s \sum\limits_{n \in \Z} \bi{q}{r+ns} + \!\!\!\!\sum\limits_{l \not\equiv r \!\! \mod s }\!\!\!\! 0.\end{align*}
Thus \begin{align*}\sum\limits_{n \in \Z}\bi{q}{r+ns}&=\frac{1}{s}\sum\limits_{j=1}^s \eta^{-2rj} \big(1+\eta^{2j}\big)^q=\frac{1}{s}\sum\limits_{j=1}^s \eta^{j(q-2r)}(\eta^{-j}+\eta^j)^q\\
&=\frac{1}{s}\sum\limits_{j=1}^s \eta^{j(q-2r)}\left(2 \cos\left(\frac{j \pi}{s}\right)\right)^q,\end{align*} 
and one takes the real part of this expression.
\end{proof}

\begin{proposition}\label{trigo} One has

 $$\dim L(\om_r)=\frac{2}{p^{k+1}} \!\!\!\sum\limits_{i=1}^{p^{k+1}-1}\!\!\!  \left(\sum\limits_{a \in A} \sin \left(\frac{i \pi (R+2a)}{p^{k+1}}\right)\right)\sin \left(\frac{i \pi}{p^{k+1}}\right) \left(2 \cos \left(\frac{i \pi}{2p^{k+1}}\right)\right)^{2m}.$$
\end{proposition}

\begin{proof}The result is easily deduced from proposition \ref{sommebinome} and lemma \ref{lemmetechnique}. With the identity $\cos A - \cos B= 2 \sin \frac{A+B}{2}\sin \frac{B-A}{2}$, one finds imediately
$$\dim L(\om_r)=\frac{1}{p^{k+1}} \!\sum\limits_{i=1}^{2p^{k+1}}\!\!  \left(\sum\limits_{a \in A} \sin \left(\frac{i \pi (R+2a)}{p^{k+1}}\right)\right)\sin \left(\frac{i \pi}{p^{k+1}}\right) \left(2 \cos \left(\frac{i \pi}{2p^{k+1}}\right)\right)^{2m}.$$

The terms corresponding to $i=p^{k+1}$ and $i=2p^{k+1}$ are equal to zero, hence we can cut the sum in two parts: one sum with $i$ between $1$ and $p^{k+1}-1$, the second sum with $i$ between $p^{k+1}+1$ and $2p^{k+1}-1$. In the second sum, set 
$i=2p^{k+1}-j$ with $j$ between $1$ and $p^{k+1}-1$, and this change of variable shows that the two sums are equal.
\end{proof}

\begin{example}\label{exemple} If $R \leq p-1$, that is $m+2-p \leq r \leq m$, we obtain
$$\dim L(\om_r)=\frac{2}{p} \sum\limits_{i=1}^{p-1} \sin \left(\frac{i \pi (m+1-r)}{p}\right)\sin \left(\frac{i \pi}{p}\right) \left(2 \cos \left(\frac{i \pi}{2p}\right)\right)^{2m}.$$

One should notice that the $\Spm$-modules $L(\om_r)$ have been investigated in \cite{Gow} for these $p-1$ particular values of $r$.

For $r=m$, the formula becomes $\dim L(\om_m)=\frac{2}{p}\sum\limits_{i=1}^{p-1} \sin^2\left(\frac{i\pi}{p}\right) \left(2+2\cos\left(\frac{i\pi}{p}\right)\right)^m\!$. We obtain $\dim L(\om_m)=2^m$ if $p=2$, $\dim L(\om_m)=\frac{1}{2}(3^m+1)$ if $p=3$ and

\noindent $\dim L(\om_m)=\frac{5-\sqrt{5}}{20}\left(\frac{5+\sqrt{5}}{2}\right)^m + \frac{5+\sqrt{5}}{20}\left(\frac{5-\sqrt{5}}{2}\right)^m + \frac{5+\sqrt{5}}{20}\left(\frac{3+\sqrt{5}}{2}\right)^m + \frac{5-\sqrt{5}}{20}\left(\frac{3-\sqrt{5}}{2}\right)^m$
if $p=5$, using $\cos\left(\frac{\pi}{5}\right)=\frac{\sqrt{5}+1}{4},\cos\left(\frac{2\pi}{5}\right)=\frac{\sqrt{5}-1}{4},\sin^2\left(\frac{\pi}{5}\right)=\frac{5-\sqrt{5}}{8},\sin^2\left(\frac{2\pi}{5}\right)=\frac{5+\sqrt{5}}{8}$.
\end{example}

In the following corollary, $d$ is a positive integer, we set $R=d+1$ and we use the notations $f$, $k$, $R_i$, $A$ as defined in definitions \ref{notation} and \ref{notationbis}. Moreover we set

$c=2^{2d} \cos^{2d} \left(\frac{\pi}{2p^{k+1}}\right) \frac{2}{p^{k+1}} \sin \left(\frac{\pi}{p^{k+1}}\right)\sum\limits_{a \in A} \sin \left(\frac{\pi (R+2a)}{p^{k+1}}\right)$.

\begin{proposition}\label{asymptotique} For given $p$ and $d$, one has  $$\dim L_{d+n}(\om_n) \mathop{\sim}\limits_{n \rightarrow \infty} c \ 4^n \cos^{2n}\left(\frac{\pi}{2p^{k+1}}\right).$$
\end{proposition}

\begin{proof}
It is a straightforward consequence of proposition \ref{trigo}.
\end{proof}

\section{Simple Weyl modules}

Now we study the Weyl modules $\Delta(\om_r)$ for $\Spm$ which are simple. In \cite{PS}, Premet and Suprunenko have obtained several necessary and sufficient conditions for a Weyl module $\Delta(\om_r)$ to be simple. We briefly recall these conditions, $v_p$ being the $p$-adic valuation on $\Q^*$.

\begin{theorem}[\cite{PS}]\label{critere} Let $1 \leq r \leq m$. The Weyl module $\Delta(\om_r)$ is simple if and only if one of the two following equivalent conditions is satisfied:
\begin{enumerate}
\item  for $1 \leq j \leq \frac{r}{2}$, one has $v_p\left(\frac{m-r+1}{j}+1\right) \leq 0$,
\item for $0 \leq j < r$ and $j \equiv r \mod 2$, $p$ does not divide $\bi{m-r+1+\frac{r-j}{2}}{\frac{r-j}{2}}$.
\end{enumerate}
\end{theorem} 

First we establish a simpler criterion for the simplicity of these Weyl modules, which will allow us to obtain all the simple Weyl modules $\Delta(\om_r)$. We use the notations of definition \ref{notation}.

\begin{theorem}\label{irred}The Weyl module $\Delta(\om_r)$ is simple if and only if  $r<2(p-R_f)p^{f}$.
\end{theorem}

\begin{proof}
  Since $L(\om_r)$ is a quotient of $\Delta(\om_r)$, the Weyl module $\Delta(\om_r)$ is simple if and only if $\dim L(\om_r)=\dim \Delta(\om_r)$. If $r<2(p-R_f)p^{f}$, from theorem \ref{theoreme} we deduce that $\dim L(\om_r)=\dim \Delta(\om_{r})$. If $r \geq 2(p-R_f)p^{f}$, set $j=r-2(p-R_f)p^{f} \geq 0$. Then $\frac{r-j}{2}=(p-R_f)p^{f} \subset m+1-j=R+2(p-R_f)p^{f}$, thus one has $\dim \Delta(\om_r) \geq \dim L(\om_{r})+\dim L(\om_{j})>\dim L(\om_{r})$ by theorem \ref{composition}.
\end{proof}

One easily sees that this criterion is equivalent to the conditions of theorem \ref{critere}.
We now generalize definition \ref{notation}.

\begin{definition}\label{def} For any integer $N \geq 2$, set $N=\sum\limits_{i=1}^h N_i p^{s_i}$ with $1 \leq N_i \leq p-1$ for any $i$ and $0 \leq s_1 < s_2 < \dots <s_h$. For $0 \leq r \leq N-1$, set $R=N-r$ and let $R=\sum\limits_{i= f}^k R_i p^i$ be the $p$-adic expansion of $R$ with $R_f \not=0$. Let $I_p(N)$ be the set of all integers $r$ satisfying $0 \leq r \leq N-1$ and $r<2(p-R_f)p^f$, and let $C_p(N)$ be the cardinal of $I_p(N)$.
\end{definition}

\begin{corollary}The Weyl module $\Delta(\om_r)$ is simple if and only if $r \in~I_p(m+~1)$.  
\end{corollary}

The following proposition gives a simple expression for $C_p(N)$. For any $x \in \R$,we note $[x]$ its integer value and $\log_p$ the logarithm to base $p$.

\begin{proposition}\label{cardinal}
 One has  $C_p(N)=(p-1)s_h+N_h=(p-1)[\log_p(N)]+\left[\frac{N}{p^{[\log_p(N)]}}\right].$
\end{proposition}

\begin{proof} Set $f(N)=(p-1)s_h+N_h$. We want to establish the equality $C_p(N)=f(N)$. Clearly one has $f(N+1)=f(N)+1$ if $N+1=dp^t$ with $1 \leq d \leq p-1$, and $f(N+1)=f(N)$ otherwise. As $C_p(2)=2=f(2)$, it is enough to prove that $C_p(N+1)=C_p(N)+1$ if $N+1=dp^t$ with $1 \leq d \leq p-1$ and $C_p(N+1)=C_p(N)$ otherwise.
One has $N=r+\sum\limits_{i=f}^k R_i p^i$, whence $N+1=(r+1)+\sum\limits_{i=f}^k R_i p^i$. Thus $r+1 \in I_p(N+1)$ implies $r \in I_p(N)$. If one has $r \in I_p(N)$, then $r+1 \in I_p(N+1)$ except when $r=2(p-R_f)p^f-1$, that is $N=2(p-R_f)p^f-1+\sum\limits_{i=f}^k R_i p^i=(p-R_f)p^f-1+p^{f+1}+R_{f+1}p^{f+1}+\dots=(p-R_f)p^f-1+p^{f+1}u$ with $u \not=0$. This happens exactly when $N+1$ has at least two non null digits in its $p$-adic expansion, and in this case the equation $r=2(p-R_f)p^f-1$ has an unique solution $r_0$ since $f$ and $R_f$ are determined by the formula $N+1=(p-R_f)p^f+p^{f+1}u$.
 Therefore if $N+1=dp^t$ with $1 \leq d \leq p-1$, one has $r' \in I_p(N+1)$ if and only if $r'=0$ or $r'-1 \in I_p(N)$, hence $C_p(N+1)=C_p(N)+1$. Otherwise,  one has $r' \in I_p(N+1)$ if and only if $r'=0$ or $r'-1 \in I_p(N) \setminus \{r_0 \}$, whence $C_p(N+1)=C_p(N)$.
\end{proof}

This proof gives an iterative construction of $I_p(N)$, but we will not use it to determine $I_p(N)$. One should note that we have obtained the formula $C_p(N+1)=C_p(N)+1$ if $N+1=dp^t$ with $1 \leq d \leq p-1$ and $C_p(N+1)=C_p(N)$ otherwise, from which we deduce that the Weyl modules for $\Spm$ with fundamental highest weight are all simple if and only if $p \geq m+1$.

The following proposition shows that some elements of $I_p(N)$ are easily obtained.

\begin{proposition}\label{sommepartielle} If $r \in I_p(N)$ and $r+p^f < N$ , then $r+p^f \in I_p(N)$. In particular one has $\sum\limits_{i=1}^j N_i p^{s_i}+dp^{s_{j+1}} \in I_p(N)$ for $0 \leq j \leq h-1$ and $0 \leq d \leq N_{j+1}-1$. 
\end{proposition}

\begin{proof}
One has $N=r+\sum\limits_{i=f}^k R_i p^i=\big(r+p^f\big)+(R_f-1)p^f+R_{f+1}p^{f+1}+\cdots$. If $R_f \geq 2$, then $r+p^f<2\big(p-(R_f-1)\big) p^f$ since $r+p^f< 2(p-R_f) p^f+p^f$, thus  $r+p^f \in I_p(N)$. Suppose $R_f=1$. As $r+p^f < N$, one has $R=p^f+\sum\limits_{i \geq g}R_i p^i$ with $g>f$ and $R_g \not=0$. Thus $r+p^f< (2p-2R_f+1) p^f < 2(p-R_g) p^g$, whence $r+p^f \in I_p(N)$.

To obtain the second part of the proposition, it is convenient to restate the first part as follows: if $N=r+R_fp^f+\sum\limits_{i>f}R_ip^i$ with $r \in I_p(N)$ and $R_f \not=0$, then $N=\big(r+p^f\big)+(R_f-1)p^f+\sum\limits_{i>f}R_ip^i$ and  $r+p^f \in I_p(N)$. Since $0 \in I_p(N)$ and $N=0+N_1p^{s_1}+N_2p^{s_2}+\cdots$, one finds successively $p^{s_1} \in I_p(N)$, $2p^{s_1} \in I_p(N)$, \dots, $N_1p^{s_1} \in I_p(N)$, $N_1p^{s_1}+p^{s_2} \in I_p(N)$,  $N_1p^{s_1}+2p^{s_2} \in I_p(N), \, \dots$. 
\end{proof}

In the following theorem, we give all the elements of $I_p(N)$. Note that the elements in the left column of the array are pairwise distincts.

\begin{theorem} Let $1 \leq r < N$. Then one has $r \in I_p(N)$ if and only if $r$ is in the following array, with $1 \leq j < h$ in each of the last four cases.
\small{$$\begin{tabular}{|l|p{3.2cm}|}
\hline
$r$ & conditions on $s$ and $d$ \\
\hline
$dp^s$ & $s < s_1$,

 $1 \leq d \leq p-1$ \\
\hline
$dp^{s_1}$ & $1 \leq d \leq N_1-1$ \\
\hline
$\sum\limits_{i=1}^j N_i p^{s_i}$ & \\
\hline
$\sum\limits_{i=1}^j N_i p^{s_i}+dp^{s_j}$ & $N_j+1 \leq d \leq p-1$ \\
\hline
$\sum\limits_{i=1}^j N_i p^{s_i}+dp^s$ & $s_j< s < s_{j+1}$,

 $1 \leq d \leq p-1$ \\
\hline
$\sum\limits_{i=1}^j N_i p^{s_i}+dp^{s_{j+1}}$ & $1 \leq d \leq N_{j+1}-1$\\
\hline
\end{tabular}$$}
\end{theorem}

\begin{proof}
We want to show that $r \geq 1$ is an element of $I_p(N)$ if and only if $r$ satisfies one of the two following conditions:
\begin{enumerate}
\item $r=dp^s$ with $s < s_1$ and $1 \leq d \leq p-1$, or $s=s_1$ and $1 \leq d \leq N_1-1$.
\item There exists $1 \leq j < h$ such that $r=\sum\limits_{i=1}^j N_i p^{s_i}+dp^s$ with $d=0$, or $s=s_j$ and $N_j+1 \leq d \leq p-1$, or $s_j< s < s_{j+1}$ and $1 \leq d \leq p-1$, or $s=s_{j+1}$ and $1 \leq d \leq N_{j+1}-1$.
\end{enumerate}

First we verify that there are exactly $(p-1)s_h+N_h-1$ integers $r \geq 1$ satisfying these conditions, as stated in proposition \ref{cardinal}. The first condition gives $(p-1)s_1+N_1-1$ values for $r$. For a given $j$, the second condition gives $1+(p-N_j-1)+(p-1)(s_{j+1}-s_j-1)+N_{j+1}-1=(p-1)(s_{j+1}-s_j)+N_{j+1}-N_j$ values. By summing we obtain the desired result.

Then it is enough to check that if $r$ satisfies one of the two preceding conditions, then one has $r \in I_p(N)$. We make the verification for the second condition, the first condition being completely similar. Set $r=\sum\limits_{i=1}^j N_i p^{s_i}+dp^s$. If $d=0$ or $s=s_{j+1}$, then we have $r \in I_p(N)$ by  proposition \ref{sommepartielle}. If $s_j< s < s_{j+1}$ and $1 \leq d \leq p-1$, or if $s=s_j$ and $N_j+1 \leq d \leq p-1$, one has $N=r+(p-d)p^s+(p-1)p^{s+1}+ \dots +(p-1)p^{s_{j+1}-1}+(N_{j+1}-1)p^{s_{j+1}}+ \cdots$, and in each case one has $r<2dp^s$.
\end{proof}

When $p=2$, we obtain the very simple statement below. 

\begin{corollary}\label{irredp2}Let $p=2$ and $m+1=\sum\limits_{i=1}^h 2^{s_i}$ with $0 \leq s_1 < s_2 < \dots <s_h$. There are exactly $s_h$ simple Weyl modules $\Delta(\om_r)$ for $\Spm$ with $1 \leq r \leq m$, and $\Delta(\om_r)$ is simple if and only if $r=2^s$ with $0 \leq s < s_1$, or $r=\sum\limits_{i=1}^j 2^{s_i}$ with $1 \leq j \leq h-1$, or $r=\sum\limits_{i=1}^j 2^{s_i}+2^s$ with $1 \leq j \leq h-1$ and $s_j < s < s_{j+1}$.
\end{corollary}

\section{Applications to modular representations of  $S_n$}

In this last section, we use the precedings results to obtain some character formulae for the symmetric group $S_n$. Note that the methods we use for the symplectic group can be applied to the spin groups, and this gives in particular a trigonometric dimension formula, completely similar to proposition \ref{trigo}, for any simple module for spin groups with highest weight a sum of two fundamental weights (see \cite{Fou}).

In \cite{Jam}, James described the composition factors of the simple $S_n$-modules associated to $p$-regular partitions of $n$ with at most two parts. We use the following notations: for any integer $i$ with $0 \leq i \leq \frac{n}{2}$, we denote by $S^{(n-i,i)}$ the Specht module associated to the partition $(n-i,i)$ and $D^{(n-i,i)}$ is the simple module associated to the same partition. 

\begin{theorem}[James]\label{James} The composition factors of $S^{(n-r,r)}$ are multiplicity free and $D^{(n-i,i)}$ is a composition factor of $S^{(n-r,r)}$ if and only if $r-i\subset n+1-2i$.
\end{theorem}

Hence we obtain the decomposition matrix of these Specht modules $S^{(n-r,r)}$ by deleting the odd lines and odd columns of the matrix $B(n+2)$: if one sets $k=r+1$ and $l=i+1$, then $D^{(n-i,i)}$ is a composition factor of $S^{(n-r,r)}$ if and only if $B(n+2)_{2k,2l}=1$. The inverse of this decomposition matrix is obtained by deleting the odd lines and odd columns of the matrix $A(n+2)$. Indeed one has $A(n+2)_{k,l}=B(n+2)_{k,l}=0$ when $k$ and $l$ have not the same parity. We obtain the following proposition, where for any $S_n$-module $M$ we denote by $[M]$ its image in the Grothendieck group $K_0(S_n)$ of finitely generated $S_n$-modules.

\begin{proposition}Let $R=n+1-2r=\sum\limits_{i=f}^kR_ip^i$ be the $p$-adic expansion of $R$ with $R_f \not=0$ and set $\delta=(p-R_f)p^ f$. One has $$\big[D^{(n-r,r)}\big]=\som{j \in J}{j \leq r}\big[S^{(n-r+j,r-j)}\big]-\som{j \in J}{j \leq r+\delta} \big[S^{(n-r+j+\delta,r-j-\delta)}\big]$$
with $J= \Big\{ j \geq 0 \ \Big\vert$ if $j=\sum\limits_{i \geq 0} j_i p^i$ is the p-adic expansion of $j$, one has $j_i=0$ when $i \leq f$ and $j_i+R_i<p$ when $i \geq f+1 \Big\}$.
\end{proposition}

\begin{proof}From the preceding discussion, one has
$$\big[D^{(n-r,r)}\big]=\sum\limits_{r-i\prec_1 n+1-2r} \big[S^{(n-i,i)}\big]-\sum\limits_{r-i\prec_{-1} n+1-2r}\big[S^{(n-i,i)}\big],$$ 
and we set $j=r-i$.
\end{proof}

This improves the last proposition of \cite{Erd} which gives the dimension of $D^{(n-r,r)}$ as a less explicit alternating sum of dimensions of Specht modules. Note that from this last proposition, one can easily deduce a trigonometric formula for $\dim D^{(n-r,r)}$ and the simple Specht modules $S^{(n-r,r)}$.


\begin{appendix}

\section{Matrices $\tilde{B}(p^n)$}

In these appendices, we establish the iterative construction of $A(p^n)$ and $B(p^n)$. Working with $\tilde{B}(p^n)$ instead of $B(p^n)$ is simpler, thus we set
 $D(p)=\Id{p}$, $F_n=
\left(\begin{array}{ccccc}
 0 & \cdots & 0 & 1 & 0 \\ 
 \vdots & \iddots &\iddots & \iddots & \vdots \\
0 & \iddots & \iddots &  & \vdots \\
1 & \iddots  & &  & \vdots \\
0 & \cdots & \cdots & \cdots & 0 
\end{array}\right) \in \mathcal{M}_n(K)$  and we define recursively the matrices $D(p^n)$ by the formula

\medskip

$D(p^{n+1})=
\left(\begin{array}{cccccc}
D(p^n) & F_{p^n}D(p^n) & 0 &\cdots &\cdots & 0 \\
0 & \ddots &\ddots & \ddots  & & \vdots \\
\vdots & \ddots  & & &  \ddots  & \vdots \\
 \vdots &  &\ddots  & & \ddots & 0 \\
\vdots & &  & \ddots & \ddots &  F_{p^n}D(p^n) \\
0 & \cdots & \cdots & \cdots & 0 & D(p^n)
\end{array}\right)$.

We want to show that $D(p^n)=\tilde{B}(p^n)$ for any $n$, what gives easily the iterative formula for matrices $B(p^n)$.
First we study the couples $(u,v)$ with $v \in \N$, $u \in \Z$, $-v \leq u \leq v$ and $\frac{v-u}{2} \subset v$. 

\begin{lemma}\label{stable} Let $u \in \Z$ and $v \in \N$ such that $-v \leq u \leq v$. One has $\frac{v-u}{2} \subset v$ if and only if $\frac{v+u}{2} \subset v$. 
\end{lemma}

\begin{proof}
If one has $a \subset b$, then clearly $b-a \subset b$.
\end{proof}

\begin{lemma}\label{addition}
Let $a,b$ be integers such that $b<p^n$ and $a \subset b$. For any $d \geq 1$ one has $a \subset b+dp^n$.
\end{lemma}

\begin{proof}
Obvious.
\end{proof}

\begin{lemma}\label{subset} Let $u, v$ be integers such that $0 \leq u \leq v$ and let $v=v_0+v_1 p + \dots + v_k p^k$ be the $p$-adic expansion of $v$ with $v_k \not=0$. One has $\frac{v-u}{2} \subset v$ if and only if $u=\pm v_0 \pm v_1 p \pm \dots \pm  v_{k-1} p^{k-1} + v_k p^k$. In particular, if $v=v_k p^k$ and $u\geq 0$, then one has $\frac{v-u}{2} \subset v$ if and only if $u=v$. 
\end{lemma}

\begin{proof}Assume $\frac{v-u}{2} \subset v$. One has a partition of $\{0, \dots, k \}$ in two sets $X$ and $Y$ such that $\frac{v-u}{2}=\sum\limits_{i \in X}v_i p^i$ and $\frac{v+u}{2}=\sum\limits_{i \in Y}v_i p^i$. As $u=\frac{v+u}{2}-\frac{v-u}{2}$, we find $u=\pm v_0 \pm v_1 p \pm \dots \pm  v_{k-1} p^{k-1} \pm v_k p^k$, and since $u \geq 0$, we must have $u=\pm v_0 \pm v_1 p \pm \dots \pm  v_{k-1} p^{k-1} + v_k p^k$. The reciprocal is clear.
\end{proof}

We show now that $D(p^n)=\tilde{B}(p^n)$ for any $n$. Since $D(p^n)$ and $\tilde{B}(p^n)$ are upper unipotent triangular matrices, it is enough to check that $D(p^n)_{u,v} =1$ if and only if $\frac{v-u}{2} \subset v$ for $1 \leq u \leq v \leq p^n$. We proceed by induction on $n$. If $n=1$ and $1 \leq u  \leq v \leq p$, then from the last part of lemma \ref{subset} one has $\frac{v-u}{2} \subset v$ if and only if $u=v$, whence $D(p)=\tilde{B}(p)$.
 
Assuming $D(p^n)=\tilde{B}(p^n)$, we determine first the couples $(u',v')$ such that $D\big(p^{n+1}\big)_{u',v'} =1$. According to the iterative construction of matrices $D(p^k)$, there are two cases to consider:
\begin{enumerate}
\item  $u'=dp^n+u$ and $v'=dp^n+v$ with $1 \leq u \leq p^n$, $1 \leq v \leq p^n$, $0 \leq d \leq p-1$ and $D(p^n)_{u,v} =1$,
\item $u'=dp^n-u$ and $v'=dp^n+v$ with $1 \leq u \leq p^n-1$,  $1 \leq v \leq p^n$, $1 \leq d \leq p-1$ and $D(p^n)_{u,v} =1$. Indeed left multiplication by $F_{p^n}$ acts on matrices almost like a horizontal symmetry. 
\end{enumerate}

Now we prove in two steps the equivalence $D\big(p^{n+1}\big)_{u',v'} =1 \Longleftrightarrow \frac{v'-u'}{2} \subset v'$.
\begin{enumerate}[(a)]
\item Suppose that $D\big(p^{n+1}\big)_{u',v'} =1$. One has $D(p^n)_{u,v} =1$ with $u,v$ defined above, thus $\frac{v-u}{2} \subset v$ by induction hypothesis. Now $\frac{v'-u'}{2}=\frac{v \pm u}{2}$ and $\frac{v \pm u}{2} \subset v$ by lemma \ref{stable}, thus $\frac{v'-u'}{2}\subset v$.  If $v < p^n$, lemma \ref{addition} shows that $\frac{v'-u'}{2} \subset v'$. If $v =p^n$, lemma  \ref{subset} gives $u=p^n$, thus we are in the first case, that is $u'=dp^n+u$ and $v'=dp^n+v$, whence $\frac{v'-u'}{2}=0 \subset v'$. Finally, we have shown that if $D\big(p^{n+1}\big)_{u',v'} =1$, then $\frac{v'-u'}{2} \subset v'$. 

\item Now we assume that $\frac{v'-u'}{2} \subset v'$ and we show that $D\big(p^{n+1}\big)_{u',v'} =1$. If $p^n \leq v'<p^{n+1}$,   we have $v'=v_0+v_1 p + \dots + v_n p^n$ and $u'=\pm v_0 \pm v_1 p \pm \dots \pm  v_{n-1} p^{n-1}+ v_n p^n$ with $v_n \not=0$ by lemma \ref{subset}. We set $v=v'-v_n p^n$, $u=\varepsilon(u'-v_n p^n) \geq 0$ with $\varepsilon \in \{1,-1\}$ and $d=v_n$. Thus one has $v'=d p^n+v$ and $u'=d p^n \pm u$ with $1 \leq d \leq p-1$, $0 \leq u \leq v  < p^n$ and $\frac{v-u}{2} \subset v$ by lemma \ref{subset}. If $1 \leq u$, then we have $D(p^n)_{u,v} =1$ by induction hypothesis, hence $D\big(p^{n+1}\big)_{u',v'} =1$.
If $u=0$, one gets readily $u'=v'=v_np^n$, thus $D\big(p^{n+1}\big)_{u',v'} =1$ since $D\big(p^{n+1}\big)$ is unipotent. If $v'<p^n$ and $\frac{v'-u'}{2} \subset v'$, one has $D(p^n)_{u',v'} =1$ by induction hypothesis and $D\big(p^{n+1}\big)_{u',v'} =1$ (case (1) with $d=0$). Finally if $v'=p^{n+1}$ and $\frac{v'-u'}{2} \subset v'$, one has $u'=v'=p^{n+1}$ from lemma \ref{subset}, whence $D\big(p^{n+1}\big)_{u',v'} =1$.
\end{enumerate}

\section{Matrices $\tilde{A}(p^n)$}

We proceed as in the preceding section. We define matrices $C(p^n)$ by setting $C(p)=\Id{p}$ and

$C\big(p^{n+1}\big)=
\left(\begin{array}{cccccc}
C(p^n) & \cdots & & & \ddots & \vdots \\
0 & \ddots & & & \ddots & C(p^n)F_{p^n}^2 \\
\vdots & \ddots & & & \ddots & -C(p^n)F_{p^n} \\
\vdots & &\ddots & & \ddots &  C(p^n)F_{p^n}^2\\
\vdots & & & \ddots & \ddots & -C(p^n)F_{p^n} \\
0 & \cdots & \cdots & \cdots & 0 &  C(p^n)
\end{array}\right)$.

We want to show that $C(p^n)=\tilde{A}(p^n)$ for any $n$, and we begin by studying the relation $\frac{v-u}{2} \prec u$.

\begin{lemma}\label{translation}Let $a \geq 0$ and $b>0$ be integers such that $a<p^n$, $b<p^n$ and $a \prec_{k} b$ with $k=\pm 1$. Let $c,d \in \N$ such that $c+d<p$. One has $a+cp^n \prec_{k} b+dp^n$.
\end{lemma}

\begin{proof}Obvious.
\end{proof}

\begin{lemma}\label{u=v} Let $1 \leq u \leq v \leq p^{n+1}$ with $\frac{v-u}{2} \prec u$. If $v=dp^n$ or $u=dp^n$ with $1 \leq d \leq p$, then $u=v$.
\end{lemma}

\begin{proof}First we show that if $v=ep^n$ with $1 \leq e \leq p$, then $u=dp^n$ with $1 \leq d \leq p$. If $u=p^{n+1}$, the claim is proved. Suppose that $u<p^{n+1}$ and let $u=u_s p^s+ \dots +u_np^n$ be the $p$-adic expansion of $u$ with $u_s \not=0$. Let $v=ep^n$ with $1 \leq e \leq p-1$ and set $z=\frac{v-u}{2}$. Since $z \prec u$, the $p$-adic expansion of $v=ep^n$ is $v=2z+u=u_s p^s+\dots$ or $v=(p-u_s)p^s+\dots$, thus $s=n$. If $v=p^{n+1}$, we obtain $s=n+1$ in the same manner.

Assume now that $u=dp^n$ with $1 \leq d \leq p-1$ and set $z=\frac{v-u}{2}$. One has $z \prec_{-1} d p^n$ or $z \prec_{1} d p^n$, and in each case this implies $z=0$ or $z \geq (p-d)p^n$. But $v=2z+u$ and 
$2(p-d)p^n+dp^n = p^{n+1}+(p-d)p^n>p^{n+1}$, thus $z=0$. If $u=p^{n+1}$, the same method gives $z=0$ or $z \geq (p-1)p^{n+1}$, and then $z=0$.
\end{proof}

\begin{lemma}\label{symetrie}
Let $u$ and $v$ be integers such that $1 \leq u <p^n$ and $u \leq v < 2p^n$, let $k=\pm 1$ and let $c \geq 1$ and $d \geq 0$ be integers with $c+d \leq p$. If $\frac{v-u}{2} \prec_k u$, then $\frac{(2cp^n-v)-u}{2} \prec_{-k} u+dp^n$.
\end{lemma}

\begin{proof}
  Set $z=\frac{v-u}{2}$ and let $u=u_s p^s+ \dots +u_{n-1}p^{n-1}$ be the $p$-adic expansion of $u$  with $u_s \not=0$.. Since $z<p^n$ and $z \prec_k u$, the $p$-adic expansion of $z$ is $z=z_s p^s+ \dots +z_{n-1}p^{n-1}$ with $z_s=0$ or $z_s=p-u_s$. One gets $\frac{(2cp^n-v)-u}{2}=cp^n-u-z=(p-u_s-z_s)p^s+\sum\limits_{i=s+1}^{n-1} (p-1-u_i-z_i)p^i+(c-1)p^n$. For any $i$ such that $s+1 \leq i \leq n-1$, one has $0 \leq p-1-u_i-z_i \leq p-1$ and $(p-1-u_i-z_i)+u_i \leq p-1$ , moreover $(c-1)+d \leq p-1$, and if $z_s=0$ (resp. $z_s=p-u_s$), then $(p-u_s-z_s)p^s=(p-u_s)p^s$ (resp. $0$). Hence $\frac{(2cp^n-v)-u}{2} \geq 0$ and  $\frac{(2cp^n-v)-u}{2} \prec_{-k} u+dp^n$.
\end{proof}

We have to show that $C(p^n)=\tilde{A}(p^n)$. The matrices $C(p^n)$ and $\tilde{A}(p^n)$ being upper unipotent triangular, it suffices to prove that for $1 \leq u \leq v \leq p^n$, one has $C(p^n)_{u,v} =1$ if and only if $\frac{v-u}{2} \prec_1 u$ and $C(p^n)_{u,v}=-1$ if and only if $\frac{v-u}{2} \prec_{-1} u$. We proceed by induction on $n$. If $n=1$, we deduce from lemma \ref{u=v} that $C(p)=\tilde{A}(p)$. We assume that $C(p^n)=\tilde{A}(p^n)$, and we determine the couples $(u',v')$ such that $C\big(p^{n+1}\big)_{u',v'} \not=0$. There are three cases: 
\begin{enumerate}
\item $u'=dp^n+u$ and $v'=dp^n+v$ with $1 \leq u \leq p^n$, $1 \leq v \leq p^n$, $0 \leq d \leq p-1$ and $C\big(p^{n+1}\big)_{u',v'}=C(p^n)_{u,v} \not=0$,
\item $u'=dp^n+u$ and $v'=(d+2+2e)p^n+v$ with $1 \leq u \leq p^n$, $1 \leq v < p^n$, $d \geq 0$, $e \geq 0$, $d+2+2e \leq p-1$ and $C\big(p^{n+1}\big)_{u',v'}=C(p^n)_{u,v} \not=0$,
\item $u'=dp^n+u$ and $v'=(d+2+2e)p^n-v$ with $1 \leq u \leq p^n$, $1 \leq v < p^n$, $d \geq 0$, $e \geq 0$, $d+2+2e \leq p$ and $C\big(p^{n+1}\big)_{u',v'}=-C(p^n)_{u,v}\not=0$.
\end{enumerate}

As before, we proceed in two steps.

\begin{enumerate}[(a)]
\item We first assume that $C\big(p^{n+1}\big)_{u',v'}=k$ with $k=\pm 1$ and $v<p^n$, and we want to prove that $\frac{v'-u'}{2} \prec_k u'$. We have $C(p^n)_{u,v}\not=0$, what implies $u \leq v <p^n$.
Suppose we are in the first two cases, that is $v'=dp^n+v$ or $v'=(d+2+2e)p^n+v$. One has $C(p^n)_{u,v}=C(p^{n+1})_{u',v'} =k$ and by induction hypothesis we obtain $\frac{v-u}{2} \prec_k u$ hence $\frac{v'-u'}{2} \prec_k u'$ by lemma \ref{translation}. In the third case, one has $C(p^n)_{u,v}=-C\big(p^{n+1}\big)_{u',v'} =-k$, hence $\frac{v-u}{2} \prec_{-k} u$ and lemma \ref{symetrie} gives $\frac{v'-u'}{2} \prec_k u'$.

Assume now that $C\big(p^{n+1}\big)_{u',v'}=k$ and $v=p^n$. We are in the first case and $u=v$ by lemma \ref{u=v}, thus $u'=dp^n+u=dp^n+v=v'$. Since $C\big(p^{n+1}\big)$ is unipotent we deduce that $C\big(p^{n+1}\big)_{u',v'}=1=k$, and we have $\frac{v'-u'}{2}=0 \prec_1 u'$.  

\item We assume that $\frac{v'-u'}{2} \prec_k u'$ and we show that $C\big(p^{n+1}\big)_{u',v'}=k$. Suppose $p^n \leq v'<p^{n+1}$ and set $z=\frac{v'-u'}{2}$. We have $z \leq v'<p^{n+1}$ and $u' \leq v'<p^{n+1}$. Let $z=z_0+ \dots +z_np^n$, $v'=v'_0+ \dots +v'_np^n$ and $u'=u'_0+ \dots +u'_n p^n$ be the $p$-adic expansions of $z,v',u'$. One has $v'=u'+2z$ and $z \prec_k u'$, and by carrying over when summing we obtain $v'_n=u'_n+2z_n$ or $v'_n=u'_n+2z_n+1$.

If $v'_n=u'_n+2z_n$, set $v=v'-v'_np^n$ and $u=u'-u'_np^n$. Suppose that $u \geq 1$, that is $u'=u'_sp^s+\dots+u'_np^n$ with $u'_s \not=0$ and $s<n$. As $z\prec_k u'$, one gets $z-z_np^n \prec_k u'-u'_np^n$, that is $\frac{v-u}{2} \prec_k u$ with $1 \leq u \leq v < p^n$, $u'=u+u'_np^n$ and $v'=v+(u'_n+2z_n)p^n$. Hence we are in the first two cases and $C(p^n)_{u,v}=k$ by induction hypothesis, whence $C\big(p^{n+1}\big)_{u',v'}=C(p^n)_{u,v}=k$. If $u=0$, we have $u'=dp^n$ with $1 \leq d \leq p-1$, and lemma \ref{u=v} gives $u'=v'$. This implies $k=1$, and $C\big(p^{n+1}\big)$ being unipotent we obtain $C\big(p^{n+1}\big)_{u',v'}=1=k$.

If $v'_n=u'_n+2z_n+1$, set $v=p^n -(v'-v'_np^n) \geq 1$ and $u=u'-u'_np^n$, and assume first that $u \geq 1$. One has $\frac{(2p^n-v)-u}{2}=z-z_np^n$ and $z-z_np^n\prec_k u'-u'_np^n=u$ as before, hence $\frac{(2p^n-v)-u}{2}\prec_k u$. One has $1 \leq u \leq 2p^n-v <2 p^n$, thus lemma \ref{symetrie} gives $\frac{v-u}{2}=\frac{(2p^n-(2p^n-v))-u}{2} \prec_{-k} u$. Finally we have $1 \leq u \leq v \leq p^n$, $C(p^n)_{u,v}=-k$, $u'=u+u'_np^n$ and $v'=(1+v'_n)p^n-v=(u'_n+2z_n+2)p^n-v$. Moreover we have $v <p^n$: otherwise we get $v=p^n$, hence $v'=v'_np^n$ and lemma~\ref{u=v} gives $u'=v'$, which contradicts $v'_n=u'_n+2z_n+1$. The conditions of the third case are satisfied, thus $C\big(p^{n+1}\big)_{u',v'}=-C(p^n)_{u,v}=-(-k)=k$. If $u=0$, lemma \ref{u=v} gives $u'=v'$, whence $k=1$ and $C\big(p^{n+1}\big)_{u',v'}=1=k$.

If $v'<p^n$, one has $C\big(p^{n+1}\big)_{u',v'}=C(p^n)_{u',v'}=k$ (case (1) with $d=0$), and if $v'=p^{n+1}$, one gets $v'=u'$ from lemma \ref{u=v} and thus $k=1$ and $C\big(p^{n+1}\big)_{u',v'}=1=k$.
\end{enumerate}

\end{appendix}

\bibliographystyle{amsplain}

\begin{thebibliography}{20}

\bibitem{AR} A. M. Adamovich and G. L. Rybnikov, \textit{Tilting modules for classical groups and Howe duality in positive characteristic}, Transform. Groups \textbf{1} (1996), 1--34.

\bibitem{BS} A. A. Baranov and I. D. Suprunenko, \textit{Branching rules for modular fundamental representations of symplectic groups}, Bull. London Math. Soc. \textbf{32} no.4 (2000), 409--420. 

\bibitem{CC} R. Carter and E. Cline, \textit{The submodule structure of Weyl modules for groups of type $A\sb{1}$}, Proceedings of the Conference on Finite Groups (Univ. Utah, 1975), Academic Press, New York (1976), 303--311.


\bibitem{Don} S. Donkin, \textit{On tilting modules for algebraic groups}, Math. Z. \textbf{212} (1993), 39--60.

\bibitem{Erd} K. Erdmann, \textit{tensors products and dimensions of simple modules for symmetric groups}, Manuscripta Math. \textbf{88} (1995), 357--386.

\bibitem{Fou} S. Foulle, \textit{Formules de caract\`eres pour des repr\'esentations irr\'eductibles des groupes classiques en \'egale caract\'eristique}, Th\`ese de doctorat, universit\'e Claude Bernard Lyon 1, juin 2004. 

\bibitem{Gow} R. Gow, \textit{Construction of $p$-1 irreducibles modules with fundamental highest weight for the symplectic group in characteristic $p$},  J. London Math. Soc. \textbf{58} no.2 (1998), 619--632.

\bibitem{Jam} G. D. James, \textit{The Representation Theory of the Symmetric Groups},  Lecture Notes in Math. \textbf{682} (1978).

\bibitem{Jan} Jantzen J., `Representations of algebraic groups', {\it Academic Press, Orlando} (1987).

\bibitem{Mat} O. Mathieu, \textit{Tilting modules and their applications}, Analysis on Homogeneous Spaces and Representation Theory of Lie Groups, Advanced Studies in Pure Mathematics \textbf{26} (2000), 145--212. 

\bibitem{PS} A. A. Premet and I. D. Suprunenko, \textit{The Weyl modules and the irreducible representations of the symplectic group with the fundamental highest weights}, Comm. Algebra  \textbf{11} (1983), 309--342.

\bibitem{Win} P. W. Winter, \textit{On the modular representation theory of the two-dimensional special linear group over an algebraically closed field}, J. London Math. Soc. ser. 2 \textbf{16} (1977), 237--252. 

\bibitem{Won} W. J. Wong, \textit{Representations of Chevalley groups in characteristic $p$}, \it Nagoya Math. J. \textbf{45} (1971), 39--78.

\end{thebibliography}

\end{document}